\documentclass[12pt]{article}

\usepackage{amssymb,amsmath,latexsym}
\voffset
-1.5cm

\usepackage{a4,amscd,graphicx,epsfig}
\usepackage[all]{xy}
\everymath{\displaystyle}

\newtheorem{teo}{Theorem}[section]
\newtheorem{lem}[teo]{Lemma}
\newtheorem{cor}[teo]{Corollary}

\newtheorem{prop}[teo]{Proposition}

\newtheorem{remark}[teo]{Remark}

\newcommand{\mr}{\mathbb{R}}
\newcommand{\mc}{\mathbb{C}}

\newcommand{\mh}{\mathbb{H}}

\newcommand{\mm}{\mathbb{M}}

\newcommand{\Cc}{{\mathcal C}}
\newcommand{\Dd}{{\mathcal D}}
\newcommand{\Ee}{{\mathcal E}}
\newcommand{\Ff}{{\mathcal F}}

\newcommand{\Ll}{{\mathcal L}}
\newcommand{\Mm}{{\mathcal M}}

\newcommand{\Pp}{{\mathcal P}}
\newcommand{\Qq}{{\mathcal Q}}

\newcommand{\Ss}{{\mathcal S}}
\newcommand{\Tt}{{\mathcal T}}

\newcommand{\sG}{\mathfrak s}
\newcommand{\lG}{\mathfrak l}

\def\ort#1{#1^\perp}

\def\Dim{\emph{Proof : }}
\def\cvd{\nopagebreak\par\rightline{$_\blacksquare$}}

\def\E#1#2{\left\langle #1,#2 \right\rangle}  
\def\fut{\mathrm{I}^+}                         
\def\pass{\mathrm{I}^-}

\def\coom{\mathrm{H}^}

\def\Ad{\mathrm{Ad}}

\def\ch{\mathrm{ch\,}}
\def\sh{\mathrm{sh\,}}

\def\d{\mathrm{d}}

\title{Linear structures on measured geodesic laminations}
\author{Francesco Bonsante$ ^{(1)}$}
\begin{document}
\maketitle
\vspace{0.7pt}\par
$^{(1)}$ Dipartimento di Matematica Applicata, via S. Maria 26, Pisa I-56127
\begin{abstract}
The space $\Mm\Ll(F)$ of measured geodesic laminations on
 a given compact closed \emph{hyperbolic} surface $F$ has a 
 canonical linear structure arising in fact from different sources in
 $2$-dimensional hyperbolic ({\it eartquake theory}) or 
 complex projective ({\it grafting}) geometry as well in 
 $(2+1)$ Lorentzian one ({\it globally hyperbolic spacetimes 
 of constant curvature}). We investigate this linear structure, by 
 showing in particular how it heavily depends on the geometric 
 structure of $F$, while to many other extents $\Mm\Ll(F)$
 only depends on the topology of $F$. This is already 
 manifest when we describe in geometric terms the sum of 
 two measured geodesic laminations in the simplest non trivial case of two
 weighted simple closed geodesics that meet each other at one point.
\end{abstract}

\section*{Introduction}
Measured geodesic laminations were pointed out by Thurston in the
late of 70's and since then they have played a fundamental role in
low-dimensional topology and geometry.\par
Given a hyperbolic surface $F$, whose topological support $S$ is a closed
orientable surface of genus $g\geq 2$, we will denote by $\Mm\Ll(F)$ the space
of measured geodesic laminations on $F$. We just mention some important
facts about that space, referring to Section~\ref{sec1} for some details.\par
- A natural \emph{action} of $\mr_{>0}$ on $\Mm\Ll(F)$ is defined by setting $t\lambda$ as the lamination
with the same support as $\lambda$ and such that the $t\lambda$-total mass of any
transverse arc is equal to the $\lambda$-mass multiplied by $t$.\par
- Every measured geodesic lamination $\lambda$ induces a positive-valued
function $\iota_\lambda$ on the space $\Cc_F$ of closed geodesics on $F$
by setting $\iota_\lambda(C)$ equal to the $\lambda$-mass of $C$. In
this way we obtain a map $\Mm\Ll(F)\rightarrow\mr_{\geq 0}^{\Cc_F}$. Such a
map is injective and we will consider on $\Mm\Ll(F)$ the topology induced by
$\mr_{\geq 0}^{\Cc_F}$.\par
- An important fact  is that there exists a topological
description of this space involving only the topology of $S$. It is possible
to define a canonical identification between $\Mm\Ll(F)$ and the space $\Mm\Ff(S)$ of
measured foliations on $S$  that in turn is homeomorphic to $\mr^{6g-6}$.\par
If $\Tt_g$ denotes the Teichm\"uller space of $S$ we can consider the set
\[
     \Tt_g\times\Mm\Ll_g=\bigcup_{[F]\in\Tt_g}\Mm\Ll(F)
\]
By previous facts it follows that $\Tt_g\times\Mm\Ll_g$ is \emph{ a trivial fiber
bundle} on $\Tt_g\times\Mm\Ll_g$ with fiber equal to $\mr^{6g-6}$.\\

Measured geodesic laminations are deeply involved in many contexts in
low-dimensional topology and geometry. In this paper we will focus on some
applications of measured geodesic laminations. We will see that in each
context we will deal with  a natural homeomorphism between $\Mm\Ll(F)$
and a $6g-6$-real vector space arises. Moreover this homeomorphism preserves
the product by positive numbers.
We will be interested in
studying the linear structures on $\Mm\Ll(F)$ obtained by such homeomorphisms.
We will see that even if they arise in different frameworks, the linear
structures induced on $\Mm\Ll(F)$ coincide (so $\Mm\Ll(F)$ is equipped with a
well-defined linear structure).

First let us introduce the  applications of measured geodesic laminations we will
deal with.
\medskip\par
1) The first one is \emph{the earthquake theory}. Given a measured geodesic lamination
$\lambda$ on a hyperbolic surface $F$ the (left or right) earthquake on $F$ along $\lambda$ is
a way to produce a new hyperbolic structure $E_\lambda(F)$ on $S$. This construction
was pointed out by Thurston \cite{Thurston:earth} and in a sense is a generalization of Dehn
twist action on $\Tt_g$. An important result due to Thurston is that given any pair of
hyperbolic structures $(F,F')$ there exists a unique (left) earthquake on $F$
relating them.
\medskip\par

2) The second application occurs in Thurston parameterization of the space of
   \emph{projective structures} on $S$. A projective structure is a maximal
   $(\mc\mathbb P^1, PSL(2,\mc))$-atlas. Thurston pointed out a geometric
   construction to associate to every hyperbolic structure $F$ equipped with a 
   measured geodesic lamination $\lambda$ a projective structure
   $Gr_\lambda(F)$ (called \emph{the grafting} of $F$ along $\lambda$) (see
   \cite{Thurston, EpMa, KulPin, McMullen} for details). Moreover the
   map
\[
   \Tt_g\times\Mm\Ll_g\ni (F,\lambda)\mapsto Gr_\lambda(F)\in \Pp(S)
\]
   turns to be a homeomorphism between $\Tt_g\times\Mm\Ll_g$ with the
   equivalence classes of (marked) projective structures.
\medskip\par
3) An important application of measured geodesic laminations occurs in
 $(2+1)$-Lorentzian geometry. Given any $\kappa\in\{0,\pm 1\}$
Mess \cite{Mess} pointed out an explicit construction to associate to every hyperbolic
surface $F$ equipped with a measured geodesic lamination $\lambda$ \emph{a maximal
spacetime} $Y_\kappa(F,\lambda)$ with \emph{constant curvature} equal to $\kappa$ and
a Cauchy surface diffeomorphic to $S$. Moreover he proved that for $\kappa\in\{0,-1\}$
his construction furnishes a parameterizations of maximal spacetimes with
constant curvature equal to $\kappa$ and Cauchy surface diffeomorphic to $S$.
An analogous statement was proved by Scannell \cite{Scannell} for the case $\kappa=1$.
Hence $\Tt_g \times \Mm\Ll_g$ arises as the fundamental structure 
encoding a priori rather different geometric objects. A clean geometric
explanation of this pervasive role of $\Tt_g \times \Mm\Ll_g$ was recently furnished
in~\cite{BenBon} by means of a general Wick rotation-rescaling theory.\\

In the present paper we focus on the fact that each of the above 
applications yields to a natural linear structure on $\Mm\Ll(F)$.
Our aim is to investigate these linear structures. In particular
we would like to give a {\it geometric description of the sum} 
of two measured geodesic laminations. We will show that it is actually
a quite difficult task. For, although we can give a description of 
 $\Mm\Ll(F)$ in purely topological terms (for instance by considering
 the atlas of $\Mm\Ll(F)$ given by train-tracks), the sum {\it heavily depends 
 on the given hyperbolic structure on $S$}. This fact already arises
 in the simplest non trivial case of two weighted simple closed geodesics
 that meet each other at one point. Also in this simplest case the
 determination of the sum lamination is not trivial at all.
 
\begin{remark}\emph{
We can consider the set of} differential quadratics \emph{$\Qq(F)$ with respect to
the conformal structure induced by the hyperbolic metric on $F$. 
Every differential quadratic induces a  horizontal
foliation on $S$, that, in turn, corresponds to a measured  geodesic
lamination on $F$. 
It is well-known that such a correspondence yields an identification of
$\Mm\Ll(F)$ with $\Qq(F)$ \cite{Masur}.
Notice that such a correspondence does not preserve the multiplication
by positive numbers (if $\lambda$ is the lamination corresponding to
$\omega$ the lamination corresponding to $t\omega$ is $t^{1/2}\lambda$).
So we will not  deal with the linear structure on $\Mm\Ll(F)$ arising from
this identification.
}\end{remark}

Let us briefly describe the contents of this paper.

In the first section we just give a brief sketch of constructions 
we have described and then we explain how it is possible to
associate to $\Mm\Ll(F)$ a liner structure.\\

In the second section we prove that linear structures corresponding to
different constructions actually coincide and in this way $\Mm\Ll(F)$ results
equipped with a canonical linear structure. Anyway let us just remark that
the topological identification between $\Mm\Ll(F)$ and $\Mm\Ll(F')$ described
by means of the canonically identification of $\Mm\Ll(F)$ (and $\Mm\Ll(F')$)
with $\Mm\Ff(S)$ \emph{ is not linear} with respect to those structures. Hence the
linear structure on $\Mm\Ll(F)$ does depend on the geometry on $F$.\\

In the third section we will deal with the problem of the sum of two measured
geodesic laminations. 
We will provide two partial results:
\medskip\par
1) We will show that the set of laminations non-intersecting a surface
with geodesic boundary $F'$ embedded in $F$  a subspace of $\Mm\Ll(F)$;
\medskip\par
2) Given two weighted simple geodesics $(C,c),(D,d)$ intersecting each other only in one point
 we will construct a sequence of weighted simple curves $$(A_n,a_n), (C_n,c_n),
 (D_n,d_n)$$ such that 
\[
    (C,c)+(D,d)=(A_n,a_n)+(C_n, c_n)+(D_n,d_n).
\]
 Moreover $A_n$ is disjoint from $C_n$ and $D_n$ and $C_n$ and $D_n$ meets
 each other in one point. The sequence is constructed by recurrence. If for
 some $n$ we have $d_n=0$ then the process ends and the sum lamination is the
 simplicial lamination given by the union of $(A_n,a_n)$ and $(C_n, c_n)$
 otherwise every term of the sum converge to a measured geodesic
 lamination. In particular $(A_n, a_n)$ tends to a weighted curve $(A_\infty,
 a_\infty)$ whereas the other terms tend to non-simplicial measured geodesic
 laminations $\lambda_\infty, \lambda'_\infty$ that are disjoint. Thus the sum
 lamination is the  union of $(A_\infty, a_\infty)$, $\lambda_\infty$ and $\lambda'_\infty$.


\section{Measured geodesic laminations}\label{sec1}
\begin{figure}[h!]
\begin{center}
\input{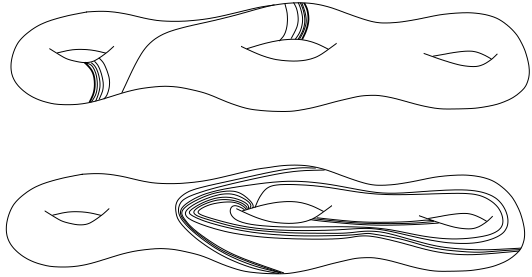}
\caption{{\small In the picture examples of non-simplicial geodesic laminations. The
    first one does not support any measure.}}\label{cap1:nonmeas:fig}
\end{center}
\end{figure}
In this paper $S$ will denote a closed orientable surface of genus $g$ and
$F$ will denote $S$ equipped with a hyperbolic metric (that is a metric of
constant curvature equal to $-1$). Moreover $\pi_1(S)$ will denote the
fundamental group of $S$ (with respect to some base point $x_0$) whereas
$\pi_1(F)$ will denote the automorphism group of a \emph{fixed} metric
covering
\[
    \mh^2\rightarrow F
\]
(so $\pi_1(F)$ is a discrete sub-group of $PSL(2,\mr)$).\\
 
A {\it geodesic lamination} $L$ on $F$ is a closed subset that is the disjoint
union of simple complete geodesics. The following list summarizes the
principal properties of geodesic laminations on a closed surface. A complete
introduction to this topic can be found in \cite{Bona, Casson}.  
\medskip\par\noindent 
1. The Lebesgue measure of a geodesic lamination is zero.  
\smallskip\par\noindent
2. If $L$ is a geodesic lamination then every point lies in a unique geodesic
contained in $L$. In particular a unique partition of $L$ in complete
geodesics exists.  
\smallskip\par\noindent 
3. The number of components of the complementary of $L$ is finite.  
\smallskip\par\noindent 
4. The number of components of $L$ is finite (but arc-connected components are
in general uncountable).  
\smallskip\par\noindent 
5. The partition of $L$ in geodesics induces a Lipschitz foliation on $L$. For this reason the
geodesics contained in $L$ are called the leaves of the lamination.
\medskip\par\noindent 
A typical example of geodesic lamination is a simple
closed geodesic or more generally a multicurve that is a disjoint union of
closed geodesics. Clearly there are more complicated geodesic laminations (see
Fig.~\ref{cap1:nonmeas:fig}).

Notice that the definition of geodesic lamination is well-founded because of
property 2. Actually in order to generalize this notion for arbitrary surfaces
it is necessary to refine the definition (see \cite{KulPin} for possible
generalizations).\\

Given a measured geodesic lamination $L$ a differentiable arc $c$ is {\it
  transverse} to $L$ if for every point $x\in L\cap c$ the leaf through $x$
is transverse to $c$. A {\it transverse measure} on $L$ is the assignment of
a Borel measure $\mu_c$ on every transverse path $c$ such that
\medskip\par\noindent 1. The support of $\mu_c$ is $L\cap c$.
\smallskip\par\noindent 2. If $c'$ is a sub-arc of $c$ then
$\mu_{c'}=\mu_c|_{c'}$.  \smallskip\par\noindent 3. If $c$ and $c'$ are
transverse paths related by an $L$-preserving homotopy then such a homotopy
sends $\mu_c$ to $\mu_{c'}$.  \medskip\par\noindent A {\it measured geodesic
  lamination} $\lambda=(L,\mu)$ is a geodesic lamination $L$ (called the
support) provided with a transverse measure $\mu$.  A simple example of
measured geodesic lamination is a {\it weighted multicurve} that is a
multicurve provided with a positive number $a(C)$ for each component $C$. If
$k$ is a transverse arc then it meets the multicurve in a finite number of
points (see Fig.~\ref{cap1:simp:fig}). The associated measure is concentrated on such points (a sum of Dirac
deltas) and the measure of any intersection point is equal to the weight of
the curve containing that point.

\begin{figure}
\begin{center}
\input{MGL_fig_cap1_simp.pstex_t}
\caption{{\small A simplicial lamination.}}\label{cap1:simp:fig}
\end{center}
\end{figure}
Carrying a transverse measure is not a property shared by all the geodesic
laminations. Actually in order to carry a transverse measure a geodesic
lamination have to satisfy certain geometric properties. For instance if $L$ is
the support of a measured geodesic lamination that it decomposes in two
sub-laminations
\[
L=L'\cup L_s
\]
such that $L_s$ is a multicurve and $L'$ does not contain any closed geodesic.
In Figure~\ref{cap1:nonmeas:fig} a geodesic lamination which does not satisfy such
a property is shown.

By multiplying a transverse measure $\mu$ by a positive number $a$ (that means
that $\mu_c$ is multiplied by $a$ for every transverse path $c$) we obtain a
new transverse measure that will be denoted by $a\mu$.  Briefly given a
measured geodesic lamination $\lambda$ and a positive number $a$ we set
$a\lambda=(L,a\mu)$.  If $\Mm\Ll(F)$ denote the set of measured geodesic
laminations then the above rule define a left action of the multiplicative
group $\mr_{>0}$ on $\Mm\Ll(F)$.

A particular lamination is the empty set. Such a lamination carries a unique
transverse measure which is the zero measure (such that the measure of any
path is zero). We will denote this degenerated measure lamination by $0$. Notice
that $0$ is the unique point fixed by $\mr_{>0}$ and the multiplication by $0$
sends every measured lamination to $0$.  
\medskip\par\noindent 
{\bf Topology  on the space of measured geodesic laminations} 
\medskip\par\noindent 
We are going to describe a suitable topology on the space $\Mm\Ll(F)$ of measured
geodesic laminations on $F$. As we are going to see this space
will be described  only in terms of topological features of $F$.

Let $\Cc$ denote the set of loops in $S$ up to free-homotopy. 
The family of closed geodesic paths,
denoted by $\Cc_F$, furnishes a complete set of representatives of the
quotient $\Cc$.
This fact will play a fundamental role in relating the geometry and the topology
of $F$. In particular it will be useful to describe  
$\Mm\Ll(F)$ just in terms of topological features of $F$. 
In what follows whenever no ambiguity arises, 
we will use $\Cc$ to indicate the
set of closed geodesics as well as the set of paths up to free homotopy.
Finally notice that the metric covering map $\mh^2\rightarrow F$ establishes a
bijection between $\Cc$ and the set of conjugacy classes of $\pi_1(F)$.

Given a geodesic lamination $L$ and a closed geodesic $C$ notice that either
$C$ is a leaf of $L$ or it is transverse to $L$.  For a given geodesic
lamination $\lambda=(L,\mu)$ let us define the intersection function
\[
\iota_\lambda:\Cc\rightarrow\mr_{\geq 0}
\]
by setting
\[
\iota_\lambda(C)=\left\{\begin{array}{ll} \mu_C(C) & \textrm{ if } C\textrm{
      is transverse to }
    L\\
    0 & \textrm{otherwise.}
                        \end{array}\right.
\]
Clearly $\iota$ is homogeneous with respect to the action of $\mr^+$, that is
\[
     \iota_{a\lambda}(C)=a\iota_\lambda(C)\qquad\textrm{ for every simple
       geodesic }C\,.
\] 
The set of simple closed geodesics of $F$ (corresponding to
the loops without self-intersections), denoted by $\Ss$, is naturally identified to the
subset of $\Mm\Ll(F)$ of curves carrying the weight $1$.  With respect to such
an identification the map $\iota_C$ associated to a simple curve $C$ is the classical
intersection form.  A classical result (\cite{Poin}) states that the intersection form
provides an embedding
\[
\Ss\rightarrow\mr_{\geq 0}^{\Cc}
\]
(actually it is possible to choose a finite number of elements of $\Cc$ in
such a way to obtain an inclusion of $\Ss$ into $\mr^N$ for $N$ sufficiently
large).  The following result is an extension of that one for general
measured geodesic laminations. In a sense it states that measured geodesic
laminations are the completion of weighted curves on $F$.
\begin{prop}\label{cap1:top_mgl:prop}
  The map 
\[
\iota:\Mm\Ll(F)\ni\lambda\mapsto\iota_\lambda\in\mr_{\geq 0}^{\Cc}
\]
is injective. Its image is the closure of the image of $\mr_+\times\Ss$  and is homeomorphic to $\mr^{6g-6}$
\end{prop}
A proof of this proposition can be found in~\cite{Penner}. 
\medskip\par\noindent
{\bf Varying the surface}
\smallskip\par\noindent
Let $F,F'$ be two hyperbolic structures on $S$  and $\iota_F,\iota_{F'}$ be the corresponding
intersection maps. Proposition~\ref{cap1:top_mgl:prop} implies that $\iota_F$
and $\iota_{F'}$ have the same image so a natural identification between
$\Mm\Ll(F)$ and $\Mm\Ll(F')$ arises by considering the map
\[
    \varphi_{FF'}=\iota_F^{-1}\circ\iota_{F'}:\Mm\Ll(F')\rightarrow\Mm\Ll(F)\,.
\]
It is possible to describe geometrically the map $\varphi_{FF'}$. Indeed given any diffeomorphism
\[
     f: F'\rightarrow F
\]
we can consider the lifting to the universal covering spaces
\[
    \tilde f: \mh^2\rightarrow \mh^2
\]
that in turn can be extended to a homeomorphism of the whole
$\overline\mh^2$ \cite{Casson}. The extension on the boundary considered up to
post-composition by elements of $PSL(2,\mr)$ does not depend on $f$ 
but only on the Teichm\"uller classes of $F$ and $F'$.

Given a lamination $L'$ of $F'$ let $\tilde L'$ denote its lifting to $\mh^2$.
For every leaf $l$ of $\tilde L'$ let $\hat f(l)$ be the geodesic with
end-points equal to the images through $\tilde f$ of the end-points of $l$.
Now the union of all $\hat f(l)$ is a geodesic lamination of $\mh^2$ invariant
under the action of $\pi_1(F)$. Thus it induces a
lamination on $F$ that we denote by $\hat f(L')$. By the above remark about
$f$ it turns out that $\hat f(L')$ does not depend on $f$ but only on $L,
F, F'$. 

Given a measured geodesic lamination $\lambda'=(L',\mu')$ on $F'$ the support
of the measured geodesic lamination $\lambda=\varphi_{FF'}(\lambda')$ is simply the
lamination $L=\hat f(L')$. In order to describe the transverse measure $\mu$
of $\lambda$ notice that it is sufficient to
describe the total mass of a  geodesic segment. Now given a geodesic segment $c$ on $F$ let
$l_-,l_+$ be the extremal leaves of $L$ cutting $c$. Let $l'_-$ and $l'_+$ be
the corresponding leaves of $L'$ and $c'$ any geodesic segment joining
them. Then we have $\mu_c(c)=\mu'_{c'}(c')$.

If $F$ and $F'$ represent the same point of the Teichm\"uller
space $\Tt_g$ then the map $\varphi_{FF'}$ is simply induced by the isometry.\\

Denote by $\Mm\Ll_g$ the image of the map $\iota$ in $\mr_{\geq 0}^{\Cc}$. As
we have seen this set depends only on $g$. Thus considering the trivial fiber bundle
\[
    \Tt_g\times\Mm\Ll_g\rightarrow \Tt_g
\]
it turns out that the fiber of a point represented by $F$ can be naturally
identified to $\Mm\Ll(F)$. Therefore $\Tt_g\times\Mm\Ll_g$ is called the fiber
bundle of measured geodesic laminations of hyperbolic surfaces of genus $g$.

\medskip\par\noindent
{\bf Intersection of measured geodesic laminations}
\smallskip\par\noindent
We have seen how it is possible to define an intersection form between a
measured geodesic lamination and a simple geodesic. Actually by using density
result of Proposition~\ref{cap1:top_mgl:prop} it is possible (see~\cite{Rees})
to define (in a
unique way) a pairing

\[
   \iota:\Mm\Ll(F)\times\Mm\Ll(F)\rightarrow\mr_+
\]
such that
\medskip\par
1. $\iota(a\lambda, a'\lambda')=a a'\iota(\lambda,\lambda')$ for every
$\lambda,\lambda'\in\Mm\Ll(F)$ and $a,a'\in\mr_{\geq 0}$;
\medskip\par
2. if $C$ and $C'$ are simple geodesics then $\iota(C,C')$ is the number of
the intersection points between $C$ and $C'$.\\

The pairing $\iota$ gives an important device to decide whether two measured
geodesic laminations transversally intersect.
\begin{teo}
Given two measured geodesic laminations $\lambda, \lambda'$ 
we have that $\iota(\lambda, \lambda')=0$ if and only if $\lambda$ and
$\lambda'$ do not transversally intersect.
\end{teo}
\cvd 
Notice that if two measured geodesic laminations $\lambda$ and $\lambda'$ do not transversally
intersect then either they are disjoint or they share some component.
Anyway the union of their supports is a geodesic lamination. 
\medskip\par\noindent
{\bf Length of a measured geodesic lamination}
\smallskip\par\noindent

Given a hyperbolic surface $F$ of genus $g$ and a closed geodesic arc $C$ we
denote by $\ell_F(C)$ its length. By identifying $\Cc$ with the set of closed
geodesics of $F$ we get a map
\[
   \ell_F:\Cc\rightarrow\mr_+
\]
called the length spectrum of $F$. It is well-known that length spectra of
hyperbolic surfaces are equal if and only if they represents the same point in
Teichm\"uller space $\Tt_g$.
Actually we can choose curves  $C_1,\ldots,C_N$ such that the map
\[
   \Tt_g \ni [F]\mapsto
   (\ell_F(c_1),\ldots\ell_F(c_N))\in\mr^N
\]
furnishes a real-analytical embedding of $\Tt_g$.\\

The length of that multicurve given
by geodesics $C_1,\ldots, C_N$ equipped with weights $a_1,\ldots, a_n$ is
simply
\[
    \ell_F(\lambda)=\sum_{i=1}^{n}a_i\ell_F(C_i)\,.
\]
\begin{prop}\label{cap1:length:prop}
There exists a unique continuous function
\[
   \ell_F:\Mm\Ll(F)\rightarrow\mr_{\geq 0}
\]
such that if $\lambda$ is a weighted multicurve then $\ell_F(\lambda)$ is its
length.
\end{prop}
See~\cite{McMullen}.

We call $\ell_F(\lambda)$ the length of the lamination $\lambda$. 


\section{Linear structure on $\Mm\Ll(F)$}
As we are going to see there are several canonical identifications of
$\Mm\Ll(F)$ with $\mr^{6g-6}$ arising from a priori very different frameworks. 
We will see that the linear structures on $\Mm\Ll(F)$ induced by
such identifications fit well so $\Mm\Ll(F)$ carries a natural linear
structure.
On the other hand we will see that the natural identification between the
spaces of measured geodesic laminations on two different hyperbolic
surfaces $F,F'$ is only a homogenous map (not linear) unless they represent
the same point of Teichm\"uller space. Thus the linear structure depends on the
geometry of $F$.
\medskip\par\noindent
{\bf Identification by flat Lorentzian geometry}
\medskip\par\noindent
Consider the isometric embedding of $\mh^2$ into the Minkowski space $\mm^3$
(that is $\mr^3$ provided with the standard scalar Minkowski scalar product
$\E{\cdot}{\cdot}$) yielded by identifying $\mh^2$ with the set
\[
    \{x| \E{x}{x}=-1\textrm{ and }x_0>0\}\,.
\]
With respect to this embedding, the isometry group of $\mh^2$ is identified to
the group of orthocronus linear transformations of $\mr^3$ preserving the
Minkowskian product.\\

\begin{figure}
\begin{center}
\input{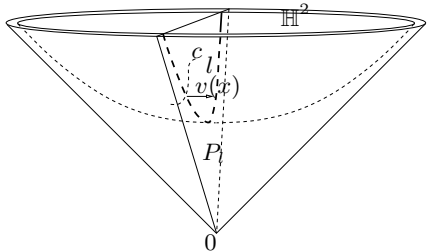}
\caption{{\small The construction of a regular domain associated to a measured
    geodesic laminations}}
\end{center}
\end{figure}
Given a measured geodesic lamination $\lambda=(L,\mu)$ on a closed hyperbolic
surface $F$ consider its lifting $\tilde\lambda=(\tilde L,\tilde\mu)$
on the universal covering $\mh^2$. Now let us fix an \emph{oriented} arc $c$ in $\mh^2$
transverse to $\tilde L$. Given a point $x\in c\cap\tilde L$ the leaf $l$
through $x$ is the intersection between a timelike plane $P_l$ and
$\mh^2$. Thus it makes sense to consider the direction orthogonal to $P_l$
that is a spacelike line. Denote by $v(x)$ the unit vector on such a line
pointing as $c$ (note that $v(x)$ depends only on $l$ and on the orientation of $c$).  
For $x$ not lying on $\tilde L$ let us put $v(x)=0$. 
In this way we have defined a function
\[
    v:c\rightarrow\mr^3
\]
that is continuous on $c\cap\tilde L$ because the foliation on $\tilde L$ is
Lipschitzian. Thus we can set
\[
     I(c)=\int_c v(x)\d\tilde\mu(x)\,.
\]
By a simple analysis of the geometry of  laminations on $\mh^2$ it is not hard
to prove that $I(c)$  depends only on the end-points of $c$ and on the
orientation of $c$. Thus given
two points $x,y$ in $\mh^2$ we choose any arc $c$ joining $x$ to $y$ and orient
it from $x$ towards $y$, and set
\[
    \rho(x,y)=I(c)\,.
\] 
Let us point some important properties of this function.
\medskip\par\noindent
1. For every $x,y\in\mh^2$ we have
\begin{equation}\label{cap2:ineq:eq}
   \E{\rho(x,y)}{y-x}\geq 0
\end{equation}
and the identity holds if and only if $x$ and $y$ lie in the same stratum of
$\tilde L$ (a stratum is either a leaf or a connected component of
$\mh^2\setminus\tilde L$). Indeed if $t$ is a point on the geodesic segment $[x,y]$ by the
choices made we have that $\E{v(t)}{y}\geq 0$ and $\E{v(t)}{x}\leq 0$
(actually the strict inequalities hold except if $v(t)=0$).
\medskip\par\noindent
2. If $\tilde L_s$ denote the lifting of the simplicial part of $L$ then we have
\begin{equation}\label{cap2:coc:eq}
    \rho(x,z)=\rho(x,y)+\rho(y,z)
\end{equation}
for every $x,z\in\mh^2$ and $y\in\mh^2-\tilde L_s$. 
\medskip\par\noindent
3. Since
$\tilde\lambda$ is invariant by the action of $\pi_1(F)$  we easily get
\begin{equation}\label{cap2:coc:eq2}
    \rho(\gamma x,\gamma y)=\gamma\rho(x,y)
\end{equation}
for every $x,y\in\mh^2$ and $\gamma\in\pi_1(F)$.\\

Fix a base point $x_0\in\mh^2-\tilde L_s$ and consider the function
\[
   \tau:\pi_1(F)\ni\gamma\mapsto\rho(x_0,\gamma x_0)\in\mr^3\,.
\]
By equations~(\ref{cap2:coc:eq}) and (\ref{cap2:coc:eq2}) we have 
\[
    \tau(\alpha\beta)=\tau(\alpha)+\alpha\tau(\beta)\,.
\]
thus $\tau$ is a cocycle of $\pi_1(F)$ taking values onto $\mr^3$ (notice that
$\mr^3$ is naturally a $\pi_1(F)$-module, since the holonomy action of
$\pi_1(F)$ extends to a linear action on $\mr^3$).

Moreover by choosing another base point $x_0'$ we obtain a new cocycle $\tau'$
that differs from $\tau$ by a coboundaries, namely
\[
   \tau'(\gamma)=\tau(\gamma)\,+\,\gamma\rho(x_0,x_0')-\rho(x_0, x_0')\,.
\]
Therefore we have defined a map
\[
    I_L:\Mm\Ll(F)\rightarrow\coom1(\pi_1(F),\mr^3)
\]
that we are going to show to be bijective.\\

Given a cocycle $\tau\in Z^1(\pi_1(F), \mr^3)$ we can associate to every
$\gamma\in\pi_1(F)$ an affine map $\gamma_\tau$ with linear part equal to $\gamma$
and translation part equal to $\tau(\gamma)$. Clearly $\gamma_\tau$ is an
isometry of the Minkowski space. Moreover the cocycle rule implies that the map
\[
   h_\tau:\pi_1(F)\ni\gamma\mapsto\gamma_\tau\mathrm{Iso}(\mm^3)
\]
is  a representation. Mess showed~\cite{Mess} that $h_\tau$ is the
holonomy of a flat spacetime homeomorphic to $F\times\mr$.

Recall that a (future complete) \emph{regular domain} is an open convex subset of
$\mr^3$ that is the intersection of the future of a non-empty family of null
planes (a null plane is a plan on which the Lorentzian form is degenerated). 
\begin{teo}\cite{Mess, Bon}\label{cap2:uniqueness:teo}
Given $\tau\in Z^1(\pi_1(F), \mr^3)$ there exists exactly one regular
domain $\Dd_\tau$ that is invariant by the action of $\pi_1(F)$ induced by
$h_\tau$. Moreover the action of $\pi_1(F)$ on $\Dd_\tau$ is free and
properly discontinuous and the quotient $Y_\tau=\Dd_\tau/\pi_1(F)$ is a
maximal globally hyperbolic spacetime homeomorphic to $F\times\mr$
\end{teo}
\cvd
We are going to sketch how it is possible to establish that the map $I_L$ is
bijective. Indeed we will show
\medskip\par\noindent
1) How to construct only in terms of $\lambda$ the regular domain $\Dd_\tau$ invariant for the cocycle
$\tau$ associated to $\lambda$.
\medskip\par\noindent
2) Given a cocycle $\tau$ how to construct a measured geodesic lamination
$\lambda$ on $F$ looking at the geometry of $\Dd_\tau$.\\

1) Let $x_0$ denote the base point of $\mh^2$ used to compute $\tau$. For $x\in\mh^2$
  let us set $u(x)=\rho(x_0,x)$. This function turns to be constant on the
  strata of $\tilde L$. Given $x\in\mh^2$ let $F(x)$ be the stratum through
  $x$ and $\partial_\infty F(x)$ the set of ideal points in the closure of
  $F(x)$ in $\overline\mh^2$. Thus we can define the set
\[
   \Omega=\bigcap_{x\in\mh^2-\tilde L}\ \bigcap_{[v]\in\partial_\infty
   F(x)}\fut(u(x)+\ort{v})\,.
\]
  By inequality~(\ref{cap2:ineq:eq}) we have that $\fut(u(x))\subset\Omega$
  and so $\Omega$ is a regular domain. On the other hand
  Equation~(\ref{cap2:coc:eq2}) implies that $\Omega$ is invariant by the
  $\pi_1(F)$-action induced by $h_\tau$. It follows that $\Omega=\Dd_\tau$.\\

2) Since $\Dd_\tau$ is a future-complete convex set it is not hard to see that
  for every point $x\in\Dd_\tau$ there exists a unique point
  $r(x)\in\partial\Dd_\tau\cap\pass(x)$ that maximizes the Lorentzian distance
  from $x$ (recall that in Minkowski space the Lorentzian distance between two
  points $x,y$ related by a timelike geodesic is simply $|x-y|=(-\E{x-y}{x-y})^{1/2}$).  
  The function $T(x)=|x-r(x)|$ carries nice properties:\\
  (i) It is the cosmological time of $\Dd_\tau$, that means that $T(x)$ is
  the $\sup$ of proper times of causal curves of $\Dd_\tau$ with
  future-endpoint equal to $x$;\\
  (ii) It is a $\mathrm C^1$-submersion and its Lorentzian gradient at $x$ is
  simply the unit timelike vector
\[
     -\frac{1}{T(x)} (x-r(x))\,;
\]
  (iii) $r(x)$ is the unique point on $\partial\Dd_\tau$ such that $x-r(x)$ is
  a spacelike support plane at $r(x)$.
  \smallskip\par\noindent
  The function $N:\Dd_\tau\rightarrow\mh^2$ given by
  $N(x)=\frac{1}{T(x)}(x-r(x))$ is called \emph{the Gauss map}. Indeed it
  coincides with the (Lorentzian) Gauss map of the level surfaces of $T$.
  The following formula can be immediately deduced by definition
\[
        x=r(x)+T(x)N(x)\,.
\]
  The image of $r$ is called \emph{the singularity} in the past of $\Dd_\tau$. By
  property (iii) it coincides with the set of points in $\partial\Dd_\tau$
  admitting a spacelike support plane. Moreover for every point
  $p$ in the singularity the set $\Ff_p=N(r^{-1}(p))\subset\mh^2$ represents
  the set of timelike directions orthogonal to some spacelike support plane at
  $p$. In particular it is not hard to see that $\Ff_p$ is a convex subset.
  Since $\Dd_\tau$ is a regular domain $\Ff_p$ turns to be the convex hull of
  its accumulation points on $\partial\mh^2$ (notice that accumulation
  points of $\Ff_p$ on $\partial\mh^2$ are the  null directions
  orthogonal to null support planes through $p$). 
Finally inequality~(\ref{cap2:ineq:eq}) implies that
  for any $p,q$ in the singularity, the geodesic of $\mh^2$ orthogonal to the
  spacelike vector $p-q$ separes $\Ff_p$ from $\Ff_q$. Thus the set
\[
    \tilde L=\bigcup_{p:\Ff_p\textrm{is a geodesic}}\Ff_p\ \cup\ \bigcup_{p:\dim\Ff_p=2}\partial\Ff_p
\]
  is a geodesic lamination of $\mh^2$. By the invariance of $\Dd_\tau$ for
  $\Gamma_\tau$ it easily follows that $\tilde L$ is invariant for
  $\Gamma$. Thus it induces a geodesic lamination $L$ on $F$\\
  In order to define a transverse measure on $\tilde L$ (invariant for $\Gamma$)
  take an arc transverse $c$ to $\tilde L$. By technical arguments \cite{BonTh}
  it is possible to prove that $u=N^{-1}(c)\cap T^{-1}(1)$ is a rectifiable
  arc. Since $r$ is Lipschitzian (with respect to the Euclidean distance) we
  can consider its derivative $\dot r$ on $u$ that turns to be a spacelike
  vector with Lorentzian length less than $1$. Thus we can set $\tilde\mu_c$ the
  direct image through $N$ of the measure $|\dot r|\d s$ where $\d s$ is the
  natural Lebesgue measure on $u$. Clearly $\tilde\lambda=(\tilde
  L,\tilde\mu)$ is a $\Gamma$-invariant measured lamination so it induces a
  geodesic lamination $\lambda$ on $F$. By construction it is not hard to see
  that the the lamination $\lambda$ induces the cocycle $\tau$.
\bigskip\par\noindent
{\bf Identification by earthquake theory}
\medskip\par\noindent

Given a measured geodesic lamination $\lambda$ on $F$ Thurston
introduced the notion of earthquake on $F$ with shearing locus $\lambda$. As
we are going to explain, it is the natural extension of the Dehn-twist action.

For the sake of simplicity we establish just the results we need, referring to the
literature on this topic to a complete introduction \cite{Thurston:earth}.
Given a weighted simple curve $(C,a)$ on a hyperbolic surface $F$  we can consider the surface $F'$
obtained by cutting $F$ along $C$ and by gluing again the geodesic boundaries
of the cutted surface after a (left) twist of length $a$. In what follows we
simply say that $F'$ is obtained by a left earthquake on $F$ with shearing
lamination $(C,a)$ and denote it by $\Ee_{(C,a)}(F)$.

Now Thurston showed that this procedure can be extended to general laminations.
\begin{teo}\cite{Thurston:earth}
There exists a continuous map
\[
    \Ee:\Tt_g\times\Mm\Ll_g\ni(F,\lambda)\mapsto\Ee_\lambda(F)\rightarrow\Tt_g
\]
such that if $\lambda$ is a weighted simple curve then $\Ee_\lambda(F)$ is the
surface describe above.

Given two elements $F,F'\in\Tt_g$ there exists a unique $\lambda\in\Mm\Ll(F)$
such that $F'=\Ee_\lambda(F)$.
\end{teo}
\cvd

Given a measured geodesic lamination $\lambda$ on a hyperbolic surface $F$
the path into the Teichm\"uller space
\[
     [0,1]\ni t\mapsto\Ee_{t\lambda}(F)
\]
turns to be differentiable. So we can associate to $\lambda$ the tangent
vector at $0$:
\[
     u_F(\lambda)=\frac{\d\,}{\d t}|_{t=0}\,\Ee_{t\lambda}(F)\,\in\mathrm T_F\Tt_g.
\]
On the other hand, since the holonomy map is a diffeomorphism of $\Tt_g$ onto
an open set of the variety of representations of $\pi_1(F)$ on $PSL(2,\mr)$ up
to conjugacy, by general facts~\cite{CaEp, Goldman} it turns out that $\mathrm T_F\Tt_g$ is
canonically identified to $\coom1_{\Ad}(\pi_1(F),\sG\lG(2,\mr))$ where $\pi_1(F)$
acts on $\sG\lG(2,\mr)$ via the adjoint representation.

In particular if $\rho_t:\Gamma=\pi_1(F)\mapsto PSL(2,\mr)$ is the holonomy
corresponding to $\Ee_{t\lambda}(F)$ then the cocycle corresponding to the
vector $u_F(\lambda)$ is simply
\[
      X_F(\lambda)[\gamma]=\frac{\d\rho_t(\gamma)}{\d t}|_{t=0}\gamma^{-1}\,.
 \]
We can  explicitly compute $X_F(\lambda)$. For every oriented geodesic $l$ of
$\mh^2$  the standard infinitesimal generator of the
group of hyperbolic transformations  with axis equal to $l$ is the element
$Y_l\in\sG\lG(2,\mr)$ such that
$\exp(Y_l)$ is the transformation with repulsive point equal to the starting
point of $l$ and translation length equal to $1$.

Now denote by $\tilde\lambda$ the lifting of $\lambda$ on $\mh^2$ and fix a
base point $x_0$. Then for $\gamma\in\Gamma$ consider the function
\[
      Y:[x_0,\gamma(x_0)]\rightarrow\sG\lG(2,\mr)
\]
such that if $t\in\tilde\lambda$ then $Y(t)$ is the standard generator of
the group of transformations with axis equal to the leaf trough $t$ (oriented as the boundary of the
half-plane containing $x_0$) and is $0$ otherwise.
Thus up to coboundaries we have~\cite{BenBon, EpMa}
\[
      X_F(\lambda)[\gamma]=\int_{[x_0,\gamma(x_0)]} Y(t)\d\tilde\mu(t)\,.
\]

Eventually we have produced a map
\[
     I_E:\Mm\Ll(F)\mapsto\in\coom1_{Ad}(\pi_1(F),\sG\lG(2,\mr))
\]
that turns to be bijective. Thus a linear structure is induced on
$\Mm\Ll(F)$. In what follows we are going to see that this structure matches
with that induced by the map $I_L$  described above.

The killing form on $\sG\lG(2,\mr)$ is a Minkowskian form so $\sG\lG(2,\mr)$
turns to be isometric to the Minkowskian space $\mr^3$. Actually there exists
a unique isometry
\[
    H: \sG\lG(2,\mr)\rightarrow\mr^3
\]
equivariant by the action of $PSL(2,\mr)$.  The map $H$ yields a isomorphism
\[
    H_*:\coom1_{Ad}(\pi_1(F),\sG\lG(2,\mr))\rightarrow\coom1(\pi_1(F),\mr^3)
\]
and we are going to see that the following diagram commutes
\[
   \begin{CD}
   \Mm\Ll(F)  @>I_E>> \coom1(\pi_1(F),\sG\lG(2,\mr))\\
   @VV Id V            @VV H_*/2 V\\
   \Mm\Ll(F) @>I_L>>  \coom1(\pi_1(F),\mr^3)\,.
   \end{CD}
\]          
Indeed if $l$ is an oriented geodesic the standard generator $X_l$
is spacelike with norm equal to $1/2$ (it is sufficient to
prove it when the axis has end-point $0,\infty$). Moreover since $H$ is
equivariant we have
\[
    \exp(X_l) H(X_l)=H(X_l)
\]
thus $H(X_l)$ is orthogonal to $l$. Finally by an explicit computation we can
see that $H(X_l)$ points outwards from the half-spaces of $\mh^2$ bounded by
$l$ and inducing the right orientation on it. 

Since $H$ is linear we have
\[
   \begin{array}{l}
   H(X_F(\lambda)[\gamma])=H\left(\int_[x_0,\gamma x_0] Y(t)\d\tilde\mu(t)\right)=\\
   \int_{[x_0,\gamma x_0]} H(Y(t))\d\tilde\mu(t)= \frac{1}{2}\tau_F(\lambda)\,.
   \end{array}
\]
{\bf Identification by using the length function}
\medskip\par\noindent
Given a hyperbolic surface $F$ and a measured geodesic lamination $\lambda$ we
have introduced the length of $\lambda$ with respect to $F$. Thus we can
consider the positive-real function $\ell$ defined on the fiber bundle of measured
laminations by setting
\[
    \ell(F,\lambda)=\ell_F(\lambda).
\]
We have:
\begin{prop}\cite{McMullen} \label{cap2:length:prop}
The map
\[
\ell:\Tt_g\times\Mm\Ll_g\rightarrow\mr_{\geq 0}
\]
is continuous. Moreover if we fix $\lambda\in\Mm\Ll_g$ then the map
\[
   u_\lambda:\Tt_g\ni F\mapsto\ell(F,\lambda)\in\mr_{\geq 0}
\]
is real-analytic.
\end{prop}
\cvd
By Proposition~\ref{cap2:length:prop} we can consider the gradient $\nabla
u_\lambda$ of $u_\lambda$ with respect to the Weil-Petersson metric of
$\Tt_g$. In this way we obtain a map
\[
   I_T: \Mm\Ll(F)\ni\lambda\mapsto\nabla u_\lambda(F)\in T_F\Tt_g\,.
\]
\begin{prop}\cite{McMullen}
Let $J$ be the endomorphism of the tangent bundle of $\Tt_g$ corresponding to
the multiplication by $i$ with respect to the complex structure of $\Tt_g$.
Then the following diagram
\[
\begin{CD}
    \Mm\Ll(F) @>I_E>> T_F\Tt_g\\
     @|                  @VV J V\\
    \Mm\Ll(F) @>I_T>> T_F\Tt_g
\end{CD}
\]
is commutative.
\end{prop}
\cvd

Thurston pointed out a construction to associate to every hyperbolic surface
$F$ equipped with a measured geodesic lamination $\lambda$ a projective
structure $Gr_\lambda(F)$ on $S$. This construction yields a parameterization
of the space of projective structures on $S$ up to projective equivalence
\[
    Gr:\Tt_g\times\Mm\Ll_g\rightarrow\Pp(S).
\] 
Given a projective structure on $S$ the maximal atlas determines a
well-defined complex structure on $S$, so we have a natural map
\[
    \Pp(S)\rightarrow\Tt_g
\]
that turns to be a holomorphic bundle.
In particular by projecting $Gr_\lambda(F)$ on $\Tt_g$ we obtain a map
\[
   \Tt_g\times\Mm\Ll_g\ni (F,\lambda)\mapsto gr_\lambda(F)\in\Tt_g\,.
\]
If we fix a pair $(F,\lambda)$ the path $t\mapsto gr_{t\lambda}(F)$ is a real
analytic path starting from $F$ so we can consider the vector

\[
     v_F(\lambda)=\frac{\d\,}{\d t}|_{t=0}\,gr_{t\lambda}(F)\,\in\mathrm T_F\Tt_g.
\]
In~\cite{McMullen} it is shown that
\[
     v_F=\nabla u_\lambda(F)
\]
so in particular we see that the grafting map induces an identification
between $\Mm\Ll(F)$ and $T_F\Tt_g$ which differ by  $I_E$ by the
multiplication by $i$ of $T_F\Tt_g$.

\section{Sum of two laminations}\label{sum-section}
We have defined a linear structure on the space $\Mm\Ll(F)$ of measured
geodesic laminations on $F$ and we have given several different
interpretations.
In this section we will take two laminations $\lambda_1,\lambda_2\in\Mm\Ll(F)$
and we will investigate what is the sum lamination
$\lambda=\lambda_1+\lambda_2$.\\
In the first part we will show that the set of measured geodesic lamination
that does not intersect an embedded surface $F'\subset F$ with geodesic
boundary is a subspace of $\Mm\Ll(F)$.\\
In the second part we give a procedure to approximate the sum lamination in
the case when the terms of the sum are simple curve meeting each other
in one point.

\medskip\par\noindent
{\bf The support of the sum lamination}
\medskip\par\noindent

Let us take $\lambda_1,\lambda_2\in\mathcal{ML}(F)$ and denote by $\lambda$ the
sum lamination $\lambda_1+\lambda_2$.
Thus the cohomological class associated with $\lambda$ is represented by the sum of
cocycles $\tau_1$ and $\tau_2$ associated to $\lambda_1$ and $\lambda_2$.\par
Let $X\subset F$ be a hyperbolic surface with totally geodesic boundary such that
the supports of $\lambda_1$ and $\lambda_2$ are contained in $X$. We will show
that the support of $\lambda$ is contained in $X$ too.\par
Let us set $F'=\overline{F-X}$ and  denote by $\tilde F'$ the inverse image
of $F'$ in $\mh^{2}$.
\begin{figure}
\begin{center}
\input{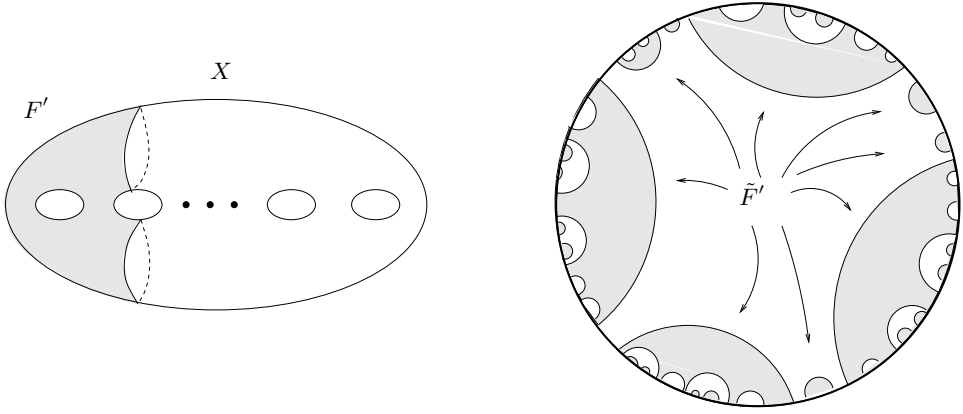}
\caption{{\small $\tilde F'$ has infinite connected components, but each
    component is open and closed in $\tilde
    F'$.}}\label{sec.5.3-suppsomma-fig}
\end{center}
\end{figure} 
Now let us fix $x_0\in\tilde F'$ and consider functions
\[
   \rho_i:\tilde F'\ni x\mapsto\int_{x_0}^x v_i(t)\mathrm d\lambda_i
   \in \mr^{2+1}\qquad\textrm{ for }i=1,2
\]
where $v_i(x)$  is the vector orthogonal to the leaf of $\lambda_i$ through $x$.
Up to adding co-boundaries we can suppose that $\tau_i(\gamma)=\rho_i(\gamma x_0)$.
Moreover since $\tilde F'$ does not intersect $\lambda_1$ and $\lambda_2$ we have that
$\rho_1$ and $\rho_2$  are locally constant functions. Since every
connected component of $\tilde F'$ is open in $\tilde F'$ it follows that
these maps are continuous.\par
Now we have to show that the function $\rho(x)=\rho_1(x)+\rho_2(x)$ has
good properties.
\begin{lem}\label{sec.5.3-support of sum-lemma}
For every $x,y\in\tilde F'$ we have that $\rho(x)-\rho(y)$ is a spacelike
vector whose dual geodesic separates $x$ from $y$. Moreover $\rho(x)-\rho(y)$  points towards $x$.
\end{lem}
\Dim
We know that $\rho_1(x)-\rho_1(y)$ and $\rho_2(x)-\rho_2(y)$ are spacelike vectors and the
corresponding dual geodesics separate $x$ from $y$. Moreover these vectors point
towards $x$. Now we have two possibilities: either the dual geodesics intersect
each other or they are disjoint. In the first  the space
generated by $\rho_1(x)-\rho_1(y)$ and $\rho_2(x)-\rho_2(y)$ is spacelike and so their sum is
spacelike.
In the second case, since they point towards $x$ we get that their scalar
product is positive. Thus it easily follows that  $\rho(x)-\rho(y)$ is
spacelike.\par
Since $\E{\rho_i(x)-\rho_i(y)}{x}\geq 0$ and $\E{\rho_i(x)-\rho_i(y)}{y}\leq
0$ the same holds for $\rho(x)-\rho(y)$. Thus if
$\rho(x)-\rho(y)\neq 0$ then  its dual geodesic separates $x$ from $y$ and
$\rho(x)-\rho(y)$ points towards $x$.
\cvd
\begin{prop}\label{sec.5.3-support of sum-prop}
Let us set $\tau=\tau_1+\tau_2$. Then we have
\[
   \mathcal D_\tau=\bigcap_{x\in\tilde F'}\fut(\rho(x)+\ort{x})\,.
\]
Moreover $\rho(x)$ lies on the singularity of $\mathcal D_\tau$.
\end{prop}
\Dim
Let $\Omega=\bigcap_{x\in\tilde F'}\fut(\rho(x)+\ort{x})$. 
First let us prove that it is a regular domain.\par
By using inequality~(\ref{cap2:ineq:eq})  we see that
$\rho(x)\in\partial\Omega$ and $\rho(x)+\ort{x}$ is a support plane through
$\rho(x)$.\par
Now let $\tilde F_1,\tilde F_2,\ldots,\tilde F_k,\ldots$ be the connected
components of $\tilde F'$ and  $\partial_\infty F_k$ be the set of ideal
points of $\tilde F_k$. Finally let us put $\rho_k=\rho(x_k)$ where $x_k\in\tilde
F_k$ (notice that $\rho_k$ does not depend on the choice of $x_k$).
It is not hard to see that
\[
  \bigcap_{x\in F_k}\fut(\rho_k+\ort{x})=\bigcap_{[v]\in\partial_\infty
  F_k}\fut(\rho_k+\ort{v})
\]
where $\partial_\infty F_k$ is the set of accumulation points of $F_k$ in the
boundary of $\mh^2$ (this equality is an easy consequence of the fact that
$F_k$ is the convex hull of $\partial_\infty F_k$).

In particular we get
\[
   \Omega=\bigcap_{k\in\mathbb N}\bigcap_{[v]\in\partial_\infty \tilde F_k}
                  \fut(\rho_k+\ort{v})\,.
\]
It follows that $\Omega$ is a regular domain and $\rho(x)$ is in the
singularity of $\Omega$. 
Moreover by construction we have that $\Omega$ is $h_\tau(\pi_1(F))$-invariant so
by Theorem \ref{cap2:uniqueness:teo} we obtain that $\Omega=\mathcal D_\tau$. 
\cvd
Let $\tilde S_1$ be the CT level surface $T^{-1}(1)$ of the domain $\mathcal
D_\tau$. Moreover let $r:\mathcal D_\tau\rightarrow\partial\mathcal D_\tau$ and
$N:\mathcal D_\tau\rightarrow\mh^{2}$ be respectively the projection on the
singularity and the Gauss map. 
By Proposition \ref{sec.5.3-support of sum-prop} we obtain the following result.
\begin{cor}\label{sec.5.3-support of sum-cor}
For every $x\in\tilde F'$ we have that $x+\rho(x)\in\tilde S_1$ and
\begin{eqnarray*}
    r(x+\rho(x))=\rho(x)\\
    N(x+\rho(x))=x\,.
\end{eqnarray*}
\end{cor}
\cvd
The lamination sum $\lambda$ is the dual lamination of the singularity
$\Sigma_\tau$ of $\Omega_\tau$.\par 
We have seen in the proof of
Proposition \ref{sec.5.3-support of sum-prop} that for every
$[v]\in\partial_\infty\tilde F _k$ the ray $\rho(x_k)+\mathbb R v$ is contained in
$\partial\mathcal D_\tau$. Thus we easily see that the plane through
$\rho(x_k)$ orthogonal to $v$ is a support plane for $\Omega_\tau$. By
definition of $\mathcal F(\rho_k)$ we obtain  $\tilde F_k\subset\mathcal
F(\rho_k)$. 
Thus $\lambda$ does not intersect the
interior of $\tilde F$.                                                                        
In particular the following corollary holds.
\begin{cor}
Let $\lambda_1,\lambda_2\in\mathcal{ML}(F)$ such that they do not intersect an
embedded surface $F'$ with totally geodesic boundary. Then the sum lamination
$\lambda=\lambda_1+\lambda_2$ does not intersect (the interior of ) $F'$.
Moreover let us fix a base point of  $x_0$ belonging to the interior of
$\tilde F'$ in the pre-image of $F'$ and
consider the cocycles $\tau,\tau_1,\tau_2$  computed with base point $x_0$
and laminations $\lambda,\lambda_1,\lambda_2$ Then we have
\[
   \tau(\gamma)=\tau_1(\gamma)+\tau_2(\gamma)\qquad\textrm{ for all }\gamma\in\pi_1(F)\,.
\]
\end{cor}
\begin{remark}\emph{
The last part of this corollary is not tautological. In fact by definition
$\tau-\tau_1-\tau_2$ is a coboundary whereas the corollary states that 
$\tau-\tau_1-\tau_2$ is zero}
\end{remark}  
\Dim
The first part of corollary is obvious. For the second one
consider the above notations. We have that
$\tau_i(\gamma)=\rho_i(\gamma(x_0))$.
On the other hand $\tau(\gamma)$ is defined by the equation
\[
     N(\gamma(x_0)+\tau(\gamma))=\gamma(x_0)\,.
\]
By Corollary \ref{sec.5.3-support of sum-cor} we have 
\[
\tau(\gamma)=\rho_1(\gamma(x_0))+\rho_2(\gamma(x_0))=
\tau_1(\gamma)+\tau_2(\gamma)\,.
\]
\cvd
\medskip\par\noindent
{\bf The sum of weighted simple curves intersecting each other only at one point}
\medskip\par\noindent

This is the simplest non trivial example of the sum of two laminations. However
we will see that even in this case the description of the sum lamination is rather involved.  
We start from simple geodesics $C$ and $D$ with weights $c$ and $d$. Then we recursively
construct a sequence of simple geodesics $A_k$, $C_k$ and $D_k$ with weights
$a_k$, $c_k$ and $d_k$ such that
\[
   (C,c)+ (D,d)=(A_k,a_k)+(C_k,c_k)+(D_k,d_k)
\]
and such that $A_k$ is disjoint from $C_k$ and $D_k$ whereas $C_k$ and $D_k$
intersect each other at one single point.
The construction ends if  $c_k$ or $d_k$ are zero for some $k$.
Otherwise the weighted curves $(A_k, a_k)$, $(C_k,c_k)$ and $(D_k,d_k)$
converge to a measured laminations $\mathcal A_\infty$, $\mathcal C_\infty$ and
$\mathcal D_\infty$. 
Moreover the transverse intersection between $\mathcal
C_\infty$ and $\mathcal D_\infty$ is zero.
Thus the union $\mathcal L_\infty=\mathcal A_\infty\cup\mathcal
C_\infty\cup\mathcal D_\infty$ is a measured lamination and we obtain that it is
the sum lamination.\par
We use the following notation: given an element $\gamma\in\pi_1(F)$ we denote
by $A_\gamma$ the axis of $\gamma$ (that is an oriented geodesic in $\mh^{2}$) and by
$C_\gamma$ the image of $A_\gamma$ in $F$.
We know that $C_\gamma$ is the unique oriented closed geodesic freely
homotopic to $\gamma$.
Finally given $\gamma,\gamma'\in\pi_1(F)$ such that $A_\gamma\cap
A_{\gamma'}\neq\varnothing$ we denote by $\theta(\gamma,\gamma')\in[0,\pi)$ the
angle between $A_\gamma$ and $A_{\gamma'}$.
\par
Now let $(C,c)$ and $(D,d)$ be two weighted simple curves
intersecting each other at one point. We have to compute 
\[
   (C,c)+(D,d)\,.
\]
Let us orient $C$ and $D$ in such a way that the angle between them is less or equal
to $\pi/2$ (if the angle is less than $\pi/2$ there are two distinct
ways to make this choice whereas if the angle is $\pi/2$ we can make every choice - i.e. there
are $4$ choices).\par
\begin{figure}[h!]
\begin{center}
\input{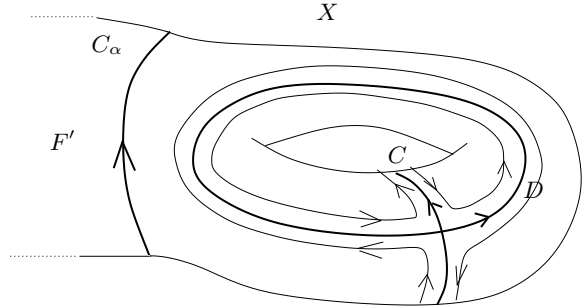}
\caption{{\small The curve $C_\alpha$ is freely homotopic to a boundary of a regular
  neighbourhood of $C_\gamma\cup C_\delta$, thus it is the boundary of a
  regular neighbourhood.}}
\end{center}
\end{figure}
Choose $\gamma,\delta\in\pi_1(F)$ such that $C_\gamma=C$ and $C_\delta=D$ as
oriented curves and $A_\gamma$ intersects $A_\delta$ at a point $p$ (we can
choose $\gamma$ arbitrarily among elements of $\pi_1(F)$ such that
$C_\gamma=C$, but the choice of $\gamma$ gives some constraints for the choice of
$\delta$).
Now let us set $\alpha=\delta^{-1}\gamma^{-1}\delta\gamma$, 
we have that $C_\alpha$ is a simple curve
which does not intersect
$C_\gamma$ and $C_\delta$. Moreover it disconnects $F$ in two regions. The
region which contains $C\cup D$ is a regular neighbourhood of this set and we
denote it by $X$. Notice that it is homeomorphic to a genus one surface with one
boundary component (i.e. a torus minus a disk). The other one, say $F'$, is a
hyperbolic surface with hyperbolic boundary.\par

Let $\tilde X$ be the component of the lifting of $X$ in $\mh^{2}$ which
contains $A_\gamma\cup A_\delta$ and $\tilde F'$ be the component of the
lifting of $F'=F-X$ which contains $A_\alpha$.
By previous paragraph we get that the support of the sum lamination
$\lambda=(C,c)+(D,d)$ is contained in $X$.
The main proposition of this section is the following one.
\begin{prop}\label{sec.5.3-sum-prop}
Suppose that $\frac{c}{d}=r(\gamma,\delta)$ where 
\[
   r(\gamma,\delta)=\frac{\cos\theta(\gamma\delta,\delta)}{\cos\theta(\gamma\delta,\gamma)}\,.
\]
Then we have
\[
  (C,c)+(D,d)=(C_\alpha, a) + (C_{\gamma\delta}, b)
\]
where $a$ and $b$ are explicit ($\mathrm C^\infty$) functions of $c,d$, the lengths of $C$ and $D$
and the angle between $C$ and $D$.
\end{prop}
The proposition is proved by a long computation. We
postpone the proof to the end of this section.
\begin{remark}\emph{
Notice that $\gamma,\delta\in\pi_1(F)$ depend (up to conjugation) on the
choice of the orientation of $C$ and $D$, in particular $\gamma\delta$ depends
on this choice. On the other hand the support of the sum lamination does not
depend on any orientation.}\par
\emph{ 
When the angle between $C$ and $D$ is less than $\pi/2$ we have two choices
for the orientation. In particular if $\gamma$ and $\delta$ 
represent $C$ and $D$ for a  given orientation then
$\gamma^{-1}$ and $\delta^{-1}$ represent $C$ and $D$ for the other one. Since
$\gamma\delta$ and $\gamma^{-1}\delta^{-1}$ are conjugated in $\pi_1(F)$ 
the result of Proposition \ref{sec.5.3-sum-prop} does not depend on  our choices.
On the other hand when the angle between $C$ and $D$ is $\pi/2$ then we can
orient geodesics so that $\gamma$ and $\delta^{-1}$ represent $C$ and $D$. 
But $\gamma\delta^{-1}$ is not conjugated to
$\gamma\delta$ nor to $(\gamma\delta)^{-1}$. However we will see that in this
case the condition is $\frac{c}{d}=1$ and the weight $b$ is equal to $0$.
}\end{remark}
\begin{figure}
\begin{center}
\input{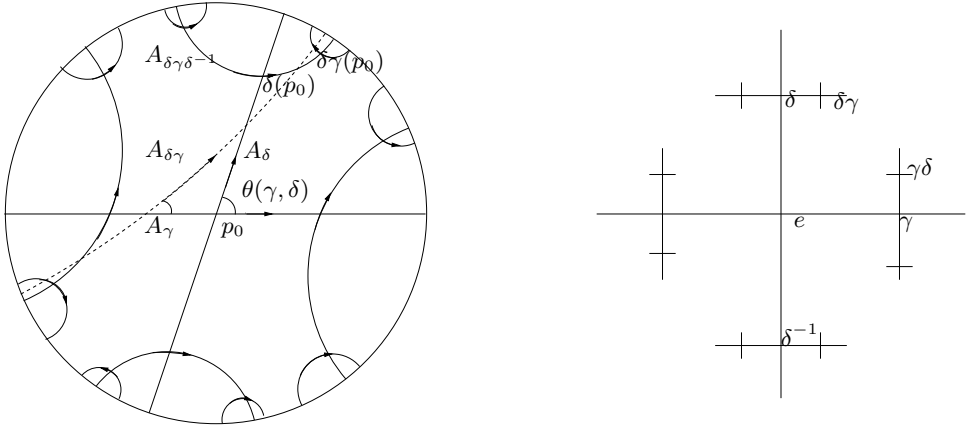}
\caption{We have  $\theta(\delta\gamma,\gamma)+\theta(\delta\gamma,\delta)<\theta(\delta,\gamma)$.}\label{sec.5.3-tree-fig}
\end{center}
\end{figure}
\begin{remark}\label{sec.5.3-angle shortens-oss}
\emph{
The stabilizer of $\tilde X$ in $\pi_1(F)$ is the free group
generated by $\gamma$ and $\delta$ (actually it is the fundamental group of
$C_\gamma\cup C_\delta$). Let us  denote this group by $\pi_1(X)$. 
Let $T$ be the component of the inverse image of $C_\gamma\cup C_\delta$
containing $A_\gamma$. It is an infinite tree such that every edge has
valence equal to $4$, see Fig.\ref{sec.5.3-tree-fig}. Vertices of $T$ are the translates of
$p_0$ by elements of $\pi_1(X)$.}\par
\emph{
Consider the Cayley graph $T'$ associated to $\pi_1(X)$: the vertices of $T'$
are the elements of $\pi_1(X)$ and two vertices are joined by an edge if they
differ by a right multiplication for $\gamma,\delta,\gamma^{-1},\delta^{-1}$.
We have that $T'$ is an infinite tree such that every vertex has valence $4$.
Moreover there exists an isomorphism of trees between $T$ and $T'$ which takes
the vertex $\eta\in\pi_1(X)$ to $\eta(p_0)$}.\par
\emph{Notice that left translations give rise to a representation of $\pi_1(X)$ into
  the group of automorphisms of $T'$. Moreover we can choose  the isomorphism
  between $T$ and $T'$ in such a way that the left multiplication corresponds to the
  natural action of $\pi_1(X)$ on $T$.}
\emph{By using this construction we can study the limit points of
  $(\delta\gamma)^n(p_0)$ for $n\in\mathbb Z$. From this analysis it follows
  that $A_{\delta\gamma}$ is like in Figure \ref{sec.5.3-tree-fig}. By looking
  at the triangle with edges on $A_\gamma\cup A_\delta\cup A_{\delta\gamma}$
  we get that
  $\theta(\delta\gamma,\gamma)+\theta(\delta\gamma,\delta)<\theta(\gamma,\delta)$ 
  (see Fig.~\ref{sec.5.3-tree-fig}).
  Thus $\theta(\delta\gamma,\gamma)\in (0,\pi/2)$ so $r(\gamma,\delta)$
  is well-defined.
}

\end{remark}
Now we will construct by recurrence sequences $\gamma_k,\delta_k\in\pi_1(F)$ and
$a_k,c_k,d_k\in\mathbb R_+$ such that
\begin{enumerate}
\item
   $(C,c)+(D,d)=(C_\alpha,a_k)+(C_{\gamma_k},c_k)+(C_{\delta_k},d_k)$.
\item
$C_\alpha$ is disjoint from $C_{\gamma_k}$ and $C_{\delta_k}$ whereas $C_{\gamma_k}$ and $C_{\delta_k}$ intersect each
other only at one point.
\item
$\alpha$ is conjugated to the commutator of $\gamma_k$ and $\delta_k$.
\item
The angle between $C_{\gamma_k}$ and $C_{\delta_k}$ is less or equal to
$\pi/2$ and $\frac{c_k}{d_k}\geq r(\gamma_k,\delta_k)$.
\item
Either there exists $k_0$ such that $d_{k_0}=0$ or the lengths of $C_{\delta_k}$
are not bounded in $\mathbb R$.
\end{enumerate}
The recurrence process ends if for some $N$ we have $d_N=0$ and in this case
we obtain that the sum $(C,c)+(D,d)$ is equal to the weighted multicurve
$(A,a_N)+(C_{\gamma_N},c_N)$.
If the process does not end then we will see that the sequence converges to
the lamination sum.\par
The first step is the following. Up to exchanging $\gamma$ with $\delta$ we
can suppose that $\frac{c}{d}>r(\gamma,\delta)$. Then let us put
\begin{eqnarray*}
\gamma_0=\gamma & \delta_0=\delta & a_0=0,\, c_0=c,\, d_0=d\,.
\end{eqnarray*}
Suppose that $\gamma_k$, $\delta_k$, $a_k$, $c_k$ and $d_k$ are defined, we have
to describe the inductive step.\par
If $d_k=0$ then we stop.
Otherwise let us consider $r_k=r(\gamma_k,\delta_k)$. We can write
\[
  (C,c)+(D,d)= (C_\alpha,a_k)+(C_{\gamma_k},c_k-r_kd_k)+ (C_{\gamma_k},r_kd_k)+(C_{\delta_k},d_k)\,.
\]
Now by applying Proposition \ref{sec.5.3-sum-prop} we get that the sum of the
two last terms is equal to 
\[
(C_\alpha,a)+(C_{\gamma_k\delta_k},b)
\]
for some $a,b\in\mathbb R_+$.
Let us put $a_{k+1}=a_k+a$. For the other curves consider the following cases.\\
If $b=0$ then put  $\gamma_{k+1}=\gamma_k$, $c_{k+1}=c_k-r_kd_k$ and
$d_{k+1}=0$.\\
If $c_k=r_kd_k$ then put  $\gamma_{k+1}=\gamma_k\delta_k$, $c_{k+1}=b$ and 
$d_{k+1}=0$.\\
If $\frac{b}{c_k-r_kd_k}\geq r(\gamma_k,\gamma_k\delta_k)$ put
\[
\left\{\begin{array}{ll}
  \gamma_{k+1}=\gamma_k\delta_k & \delta_{k+1}=\gamma_k\\
  c_{k+1}=b & d_{k+1}=c_k-r_kd_k.
\end{array}\right.\,.
\]
Finally if $\frac{b}{c_k-r_kd_k}<r(\gamma_k,\gamma_k\delta_k)$ put
\[
\left\{\begin{array}{ll}
  \gamma_{k+1}=\gamma_k & \delta_{k+1}=\gamma_k\delta_k\\
  c_{k+1}=c_k-r_kd_k & d_{k+1}=b
\end{array}\right.\,.
\]
Since $C_{\gamma_k}\cap C_{\delta_k}$ is a single  point  the 
same happens for $C_{\gamma_{k+1}}\cap C_{\delta_{k+1}}$.
Moreover by Remark \ref{sec.5.3-angle shortens-oss} the angle between
$A_{\gamma_{k+1}}$ and $A_{\delta_{k+1}}$ is smaller than the angle between
$A_{\gamma_k}$ and $A_{\delta_k}$.
Finally notice that the commutator of $\gamma_{k+1}$ and $\delta_{k+1}$ is
conjugated to $\alpha$. Thus $C_\alpha$ is disjoint from $C_{\gamma_{k+1}}$
and $C_{\delta_{k+1}}$.\par
By using these facts we can see that this sequence verifies properties 1-4.
Suppose that the sequence is infinite. 
Since $\delta_k$'s are all different, they form a divergent sequence in $\pi_1(F)$.
On the other hand since they are word in
$\gamma$ and $\delta$ with all positive exponents we get that 
$A_{\delta_k}$ have endpoints in the opposite segments of
$\partial\mh^{2}-(A_\gamma\cup A_\delta)$. Thus the translation length of $\delta_k$ 
goes to infinity.
\begin{lem}\label{sec.5.3-compactness-lemma}
Suppose that the sequence $\{\gamma_k,\delta_k,a_k,c_k,d_k\}$ is infinite. Let
us take $p_0\in\tilde F'$ and  $\beta\in\pi_1(F)$. Suppose that the geodesic segment 
$[p_0,\beta p_0]\subset\mh^{2}$ is not contained in the axis of any
element of $\pi_1(F)$ in the conjugacy class of  $\gamma_k$ and $\delta_k$. 
Let $N_k$ (resp. $M_k$) be the cardinality of the
intersection of $[p_0,\beta p_0]$ with $\tilde C_k$ (resp. $\tilde D_k$)
where $\tilde C_k$ (resp. $\tilde D_k$) is the lifting in $\mh^{2}$ of the
curve $C_k$ (resp. $D_k$). Then 
there exists $C\in\mathbb R_+$ such that $N_kc_k\leq C$ and $M_k d_k\leq C$.\par
Moreover  $a_k$'s are bounded.
\end{lem}
\begin{figure}
\begin{center}
\input{MGL_fig_cap3_triangle.pstex_t}
\caption{{\small The angle at $q$ of the triangle $pqp'$ is equal to
    $\cos^{-1}(\E{u_i}{w_j})$.}}\label{sec.5.3-triangle-fig}
\end{center}
\end{figure}
 \Dim
The cocycle associated to the sum lamination $(C,c)+(D,d)$
computed with starting point $p_0$ is equal to the cocyle associated to
$(C_\alpha,a_k)+(C_{\gamma_k},c_k)+(C_{\delta_k},d_k)$ computed with starting point $p_0$.
Let $\tau$ be such a cocycle, we know that
\[
  \tau(\beta)=a_k \sum_{i=1}^K v_i +c_k\sum_{i=1}^{N_k} w_i +
  d_k\sum_{i=1}^{M_k} u_i
\]
where $K,N_k,M_k$ are respectively the cardinalities of the intersection of
$[p_0,\beta p_0]$ with
$\tilde C_\alpha$,$\tilde C_k$ and $\tilde D_k$, whereas $v_i$, $w_i$ and
$u_i$ are respectively the unit vectors orthogonal to $\tilde C_\alpha$, $\tilde
C_k$ and $\tilde D_k$  pointing towards $\beta p_0$.
The geodesic corresponding to $v_i$ is disjoint from all the geodesics corresponding
to $v_j$, $w_j$, $u_j$. 
By an usual argument  we get $\E{v_i}{v_j}\geq 1$, $\E{v_i}{w_j}\geq
1$ and $\E{v_i}{u_j}>1$.
In the same way we have that $\E{w_i}{w_j}\geq 1$ and $\E{u_i}{u_j}\geq 1$.
Now we claim that there exists a number $L$ (independent of
$n$) such that the number of couples $(u_i,w_j)$ such that $\E{u_i}{w_j}<0$ is
less than $L$.
By the claim we get that
\[
  \E{\tau(\beta)}{\tau(\beta)}\geq (k a_k)^2 + (N_k c_k)^2 + (M_k d_k)^2 -
  2L c_kd_k
\]
(indeed if $\E{u_i}{w_i}<0$ then by construction the dual geodesics intersect
each other and so $\E{u_i}{w_i}\geq-1$). Thus the lemma follows from the claim.\par
Let us prove the claim. Suppose that $\E{u_i}{w_j}<0$, then 
the corresponding geodesics intersect each other at a point $q$. On the other hand let
$p\in\mh^{2}$ ($p'\in\mh^{2}$) be the intersection of the segment $[p_0,\beta p_0]$ with the
geodesic corresponding to $u_i$ (resp. $w_i$). Since $\E{u_i}{w_j}<0$  the
angle at $q$ of the hyperbolic
triangle  $qpp'$ is greater than $\pi/2$ (see Fig.~\ref{sec.5.3-triangle-fig}). So the distance between
$q$ and the segment $[p_0,\beta p_0]$ is less than the length of the
segment. Let $H$ be the set of points whose distance from $[p_0,\beta p_0]$ is
less than the length of this segment. We have that $H$ is a compact set so
that it intersects just a finite number $L$ of the translates of a fixed fundamental domain for
the action of $\pi_1(F)$.\par
We will see that $L$ works. In fact we have that the point $q$
projects on the intersection of $C_{\gamma_k}$ and $C_{\delta_k}$. Thus $q$
runs over a set of  $L$ elements of $\mh^{2}$. On the other hand if we
choose $q$ in this set the lifting of $C_{\gamma_k}$ (and $C_{\delta_k}$) 
passing through $q$ is unique, so the couple $(u_i,w_j)$ is determined by
$q$.
\cvd     
By Lemma \ref{sec.5.3-compactness-lemma} it follows that the families of
weighted multi-curves $\{(C_{\gamma_k}, c_k)\}$ and $\{(C_{\delta_k},d_k)\}$ are
relatively compact in $\mathcal{ML}(F)$. Thus up to passing to a subsequence
we can suppose that they respectively converge to two measured laminations
$\lambda'_\infty$ and $\lambda''_\infty$ and moreover $a_k\rightarrow
a_\infty$.

\begin{prop}
We have 
\[
  (C,c)+(D,d)=(C_\alpha,a_\infty) + \lambda'_\infty +\lambda''_\infty\,.
\]
Moreover we have
\begin{eqnarray*}
\iota(C_\alpha,\lambda'_\infty)=0 & \iota(C_\alpha,\lambda''_\infty)=0 & 
\iota(\lambda'_\infty,\lambda''_\infty)=0\,.
\end{eqnarray*} 
where $\iota:\Mm\Ll(F)\times\Mm\Ll(F)\rightarrow\mr_+$ is the intersection pairing.
\end{prop}
\Dim
The first statement follows from the construction of the sequence.
The intersection of $C_\alpha$ and $\lambda'_\infty$
(resp. $\lambda''_\infty$) is zero because of
we have that $\iota(C_\alpha,(C_{\gamma_k},c_k))=0$
(resp. $\iota(C_\alpha,(C_{\delta_k},d_k))=0$)
 and the intersection is
continuous function of $\mathcal{ML}(F)\times\mathcal{ML}(F)$.
Finally notice that
\[
 \iota(\lambda'_\infty,\lambda''_\infty)=\lim_{k\rightarrow\infty}\iota((C_{\gamma_k},c_k),(C_{\delta_k},d_k))=
\lim_{k\rightarrow\infty}c_kd_k\,.
\]
We have noticed that the length of $C_{\delta_k}$ goes to infinite so
 $d_k$ goes to zero. On the other hand we know that $c_k$ is bounded in
$\mathbb R$ so the proof is complete.
\cvd   
Since the geometric intersection between $\lambda'_\infty$ and
$\lambda''_\infty$ is zero we see that their supports have empty transverse
intersection. So the union of their supports is a geodesic lamination
too. Thus this lamination can be endowed with a transverse measure so that the
corresponding measure geodesic lamination $\lambda_\infty$ is equal to
$\lambda'_\infty+\lambda''_\infty$.
Since $\lambda_\infty$ is disjoint from $C_\alpha$ it follows that the union
of these measured laminations gives the sum lamination.
\begin{remark}\emph{
Notice that the sum lamination  has always simplicial components. Since
the sequence $a_k$ is increasing we have $a_\infty\neq 0$ so that
$(C_\alpha, a_\infty)$ is a simplicial sub-lamination of the sum lamination.
}\end{remark}
In the last part of this section we prove Proposition \ref{sec.5.3-sum-prop}.
Given a hyperbolic transformation $\alpha\in\mathrm{SO}(2,1)$ we denote by
$x^0(\alpha)\in\mr^{2+1}$ the unit spacelike vector of $\mr^{2+1}$ corresponding to
$A_\alpha$,  such that
it induces on $A_\alpha$ the orientation from the repulsive fixed point towards
the attractive fixed point.
The following lemma is a technical result which we need for the proof of
Proposition \ref{sec.5.3-sum-prop}.
\begin{lem}\label{sec.5.3-comp.-lemma}
Let $\gamma,\delta\in\pi_1(F)$ be  such that $C_\gamma$ and $C_\delta$ are two
simple curves which intersect each other at one single point.
Let $\alpha=\delta^{-1}\gamma^{-1}\delta\gamma$ and let $W$ be the subspace of
$\mr^{2+1}$ generated by $x^0(\delta \gamma), x^0(\gamma)$ and
$x^0(\alpha)-\delta x^0(\alpha)$. Then the dimension of $W$ is $2$.
\end{lem}
\Dim
Consider matrices
\begin{eqnarray*}
M(l)=\left(\begin{array}{ccc}
           \ch l  & \sh l & 0\\
           \sh l  & \ch l & 0\\
           0      & 0     & 1
           \end{array}\right) \\
R_\theta=\left(\begin{array}{ccc}
            1    & 0           & 0          \\
            0    & \cos\theta  & -\sin\theta\\
            0    & \sin\theta  & \cos\theta
            \end{array}\right)\,.
\end{eqnarray*}
We can choose coordinates in such a way  that $\gamma=M(l)$ and $\delta=R_\theta M(m)
R_{-\theta}$ where $l$ (resp. $m$) is the  length of $C_\gamma$ (resp. $C_\delta$)
and $\theta$ is the angle between $C_\gamma$ and $C_\delta$.
Thus we have that
\begin{eqnarray*}
x^0(\gamma)=\left(\begin{array}{l} 0\\0\\1\end{array}\right) &   
x^0(\delta)=\left(\begin{array}{l} 0\\-\sin\theta\\
    \cos\theta\end{array}\right)\,.
\end{eqnarray*}
By an explicit computation we have that
\[
 w=\left(\begin{array}{l}
         \sin\theta(\ch m-1)\sh l\\
         -\sin\theta(\ch m -1)(\ch l+1)\\
         \sh l\sh m + \cos\theta(\ch l+1)(\ch m-1)
         \end{array}\right)
\]
is a  generator of $\ker(\delta\gamma-1)$. 
In order to compute a generator of $\ker(\alpha-1)$ notice that
$\ker(\alpha-1)=\ker (\delta\gamma-\gamma\delta)$ The latter is a
skew-symmetric 
matrix so it is straightforward to compute its kernel.
By performing such a  computation it turns out that $\ker(\alpha-1)$ is
generated by
\[
v=\left(\begin{array}{l}
        \sh l\sh m + (\ch l-1)(\ch m-1)\cos\theta\\
        -\sh m(\ch l-1) - \sh l(\ch m-1)\cos\theta\\
         -\sin\theta\sh l(\ch m-1)
        \end{array}\right)\,.
\]
By an explicit computation we have 
\[
\delta v=  \left(\begin{array}{l}
        \sh l\sh m  - (\ch l-1)(\ch m-1)\cos\theta\\
        -\sh m(\ch l-1) + \sh l(\ch m-1)\cos\theta\\
         +\sin\theta\sh l(\ch m-1)
        \end{array}\right)\,.
\]
So we obtain
\[
v-\delta v =2(\ch m-1)\left(\begin{array}{l}
                            (\ch l-1)\cos\theta\\
                            -\sh l\cos\theta\\
                            -\sh l\sin\theta\end{array}\right)\,.
\]
Notice that $W$ is generated by $x^0(\gamma),w,v-\delta v$.
On the other hand an easy computation shows 
\[
\mathrm{det}\left[\begin{array}{lll}
                  0     &    \sin\theta(\ch m-1)\sh l   & (\ch
                  l-1)\cos\theta\\
                  0     & -\sin\theta(\ch m-1)(\ch l+1) & -\sh l\cos\theta\\
                  1     & \sh l\sh m +\cos\theta(\ch l+1)(\ch m-1) & \sh
                  l\sin\theta
                   \end{array}\right]=0\,.
\]
\cvd
\emph{Proof of Proposition \ref{sec.5.3-sum-prop}:}
We use the notation introduced above. In particular let $p_0\in\tilde F'$ be a
base point. For a given weighted curve $(A,a)$ we will denote by
$(A,a)[\gamma]\in\mr^{2+1}$ the value at $\gamma$ of the cocycle
corresponding to $(A,a)$ computed with base point $p_0$.
\begin{figure}
\begin{center}
\input{MGL_fig_cap3_cocycle.pstex_t}
\caption{}\label{sec.5.3-cocycle-fig}
\end{center}
\end{figure}
Now we want to show that under the assumptions of the proposition there exist
\emph{positive} constants $a,b$ such that

\begin{equation}\label{sec.5.3-target-eq}
   (C,c)[\beta] + (D,d)[\beta] = (C_\alpha,a)[\beta] +
   (C_{\delta\gamma},b)[\beta] \qquad\textrm{ for all }\beta\in\pi_1(F)\,.
\end{equation}
By an application of Van Kampen theorem we know that  $\pi_1(F)$
is the amalgamation of the stabilizer $\pi_1(F')$ of $\tilde F'$ with the
stabilizer of $\tilde X$ along the stabilizer of the geodesic $\tilde F'\cap
\tilde X$.
We have that the stabilizer of $X$ is the free group on $\gamma$ and
$\delta$ whereas the stabilizer of $\tilde F'\cap X$ is the group generated by
$\alpha$.\par
Notice that for all $\beta\in\pi_1(F')$ all terms
involved in expression (\ref{sec.5.3-target-eq}) are zero. Thus it is sufficient
to find $a,b\in\mathbb R_+$ such that
\begin{equation}\label{sec.5.3-target2-eq}
\left\{\begin{array}{l}
  (C,c)[\gamma]+(D,d)[\gamma]=(C_\alpha,a)[\gamma]+(C_{\delta\gamma},b)[\gamma]\\
  (C,c)[\delta]+(D,d)[\delta]=(C_\alpha,a)[\delta]+(C_{\delta\gamma},b)[\delta]\,.
\end{array}\right.
\end{equation}
Thus let us compute the terms in this expression.
By an analysis of Fig.~\ref{sec.5.3-cocycle-fig} we obtain
\[
\left\{\begin{array}{ll}
  (C,c)[\gamma]=0 & (D,d)[\gamma]=-dx^0(\delta) \\
  (C_\alpha,a)[\gamma]=a(1-\gamma)x^0(\alpha) &
  (C_{\delta\gamma},b)[\gamma]=-bx^0(\gamma\delta)=-b\gamma( x^0(\delta\gamma))\\
  (C,c)[\delta]=x^0(\gamma) & (D,d)[\delta]=0 \\
  (C_\alpha,a)[\delta]=a(1-\delta)x^0(\alpha) &
  (C_{\delta\gamma},b)[\delta]=b(x^0(\delta\gamma))\,.
\end{array}\right.
\]
Thus equation (\ref{sec.5.3-target2-eq}) is equivalent to the system
\begin{equation}\label{sec.5.3-target-eq2}
\left\{\begin{array}{l}
  a(1-\delta)x^0(\alpha)+b x^0(\delta\gamma) = c x^0(\gamma)\nonumber\\
  a(1-\gamma)x^0(\alpha)-b\gamma x^0(\delta\gamma)=-d x^0(\delta)\,.
\end{array}\right.
\end{equation}
By Lemma \ref{sec.5.3-comp.-lemma}  each equation of this system 
has a unique solution depending linearly on the weights $c$ and $d$.
Thus there exists a real number $k$ such that the solution of the first
equation coincides with the solution of the second one (i.e. the system
(\ref{sec.5.3-target-eq2}) has solution) if and only if $c/d=k$.
\par
In order to compute the coefficient $k$ notice that it is sufficient to
compute $b$ in both the equations. Now take the first equation and consider
the scalar product of each terms with $x^0(\delta)$. We have
\[
   b\E{x^0(\delta\gamma)}{x^0(\delta)}=c\E{x^0(\gamma)}{x^0(\delta)}
\]
so 
\[
   b=c\frac{\E{x^0(\gamma)}{x^0(\delta)}}{\E{x^0(\delta\gamma}{x^0(\delta)}}=
   c\frac{\cos\theta(\gamma,\delta)}{\cos\theta(\delta,\gamma\delta)}\,.
\]
On the other hand by taking the scalar product of the second equation with
$x^0(\gamma)$ we get
\[
   -b\E{\gamma x^0(\delta\gamma)}{x^0(\gamma)}=-d\E{x^0(\delta)}{x^0(\gamma)}
\]
so that we have 
\[
b=d\frac{\E{x^0(\gamma)}{x^0(\delta)}}{\E{x^0(\gamma)}{x^0(\gamma\delta)}}
=d\frac{\cos\theta(\gamma,\delta)}{\cos\theta(\gamma,\gamma\delta)}\,.
\]
Thus the system (\ref{sec.5.3-target-eq2}) has a solution if and only if
\[
    \frac{c}{d}=\frac{\cos\theta(\delta,\delta\gamma)}{\cos\theta(\gamma,\delta\gamma)}\,.
\]
Notice that we can argue this result in the case
$\theta(\gamma,\delta)\neq\frac{\pi}{2}$.
On the other hand since $k$ depends continuously on $\theta(\gamma,\delta)$
we have that this is true also in the case $\theta(\gamma,\delta)=\frac{\pi}{2}$.\par
Now we have to show that in the case $\frac{c}{d}=k$ the solutions $a,b$ of
equations
(\ref{sec.5.3-target-eq2}) are non-negative. From the above
calculation it follows that $b\geq 0$ and $b=0$ if and only if
$\theta(\gamma,\delta)=\pi/2$.
In order to compute $a$ notice that the second equation in
(\ref{sec.5.3-target-eq2}) is equivalent to the following
\[
  a(\delta-\delta\gamma)x^0(\alpha)-bx^0(\delta\gamma)=-dx^0(\delta)\,.
\]
By summing this equation to the first one of (\ref{sec.5.3-target-eq2}) we get
\[
  a(1-\delta\gamma)x^0(\alpha)=c x^0(\gamma)-d x^0(\delta)\,.
\]
(notice that by taking the scalar product of this equation with
$x^0(\delta\gamma)$ we recover the condition on $k$).
Thus by taking the scalar product with $x^0(\alpha)$ we obtain
\[
   a(1-\E{\delta\gamma x^0(\alpha)}{x^0(\alpha)})= c
   \E{x^0(\gamma)}{x^0(\alpha)}-d\E{x^0(\delta)}{x^0(\alpha)}\,.
\]
Now a careful analysis of Figure \ref{sec.5.3-cocycle-fig} shows that
\[
\begin{array}{l}
\E{\delta\gamma x^0(\alpha)}{x^0(\alpha)}=\E{\gamma x^0(\alpha)}{\delta^{-1}
  x^0(\alpha)}<0\\
\E{x^0(\alpha)}{x^0(\gamma)}>0\\ 
\E{x^0(\alpha)}{x^0(\delta)}<0\,.
\end{array}
\]
Thus it follows that $a >0$.
\cvd
\begin{remark}
\emph{
If $\theta(\gamma,\delta)=\frac{\pi}{2}$ the process ends at first step. 
Thus it turns out that if the angle between the geodesics $C_\gamma$ and
$C_\delta$ is $\pi/2$ then the sum is always a weighted multicurve (actually it
has either one component $(C_\alpha,a)$ or two components
$(C_\alpha,a)+(C_\gamma,c-kd)$).}
\end{remark}


\begin{thebibliography}{99}
\bibitem{BenBon} R.~Benedetti and F.~Bonsante,~\emph{Wick rotations in
    3D-gravity: $\Mm\Ll(\mh^2)$-spacetimes}. Preprint.
\bibitem{Bona}F.~Bonahon,~\emph{Geodesic laminations on
    surfaces}. Contemp. Math., \textbf{269} (1997), 1--37.
\bibitem{Bon}F.Bonsante,~\emph{Flat Spacetimes with Compact Hyperbolic Cauchy
    Surfaces}. To appear.
\bibitem{BonTh}F.~Bonsante,~\emph{Deforming the Minkowskian cone of a closed
    hyperbolic manifold} Ph.D. Thesis, Pisa, 2005.
\bibitem{Casson} A.~Casson,~\emph{Automorphisms on surfaces after Nielsen and
    Thurston}. London Mathematical Society Student Texts, 9. Cambridge
    University Press, 1988.
\bibitem{CaEp}R.~D.~Canary and D.~B.~Epstein,~\emph{Notes on notes of Thuston}.In
    \emph{Analytical and geometric aspects of hyperbolic space},
    3--92. Cambridge University press, 1987.
\bibitem{EpMa} D.~B.~A.~Epstein and A.~Marden,~\emph{Convex hulls in hyperbolic
    space, a theorem of Sullivan, and measured pleated surfaces}. In
    \emph{Analytical and geometric aspects of hyperbolic space},
    113--254. Cambridge University press, 1987.
\bibitem{Goldman} W.~Goldman,~\emph{The simpletic nature of fundamental groups
    of surfaces}. Adv.~in~Math., \textbf{54} (1984),200-225.
\bibitem{Masur} J.~Hubbard and H.~Masur,~\emph{Quadratic differentials and
    foliations}. Acta Math., \textbf{142} (1979), 221--274.
\bibitem{KulPin} R.~Kulkarni and U.~Pinkall,~\emph{A canonical metric for
    Moebius structures and its applications}. Math.~Z., \textbf{216} (1984), 89--129.
\bibitem{McMullen}C.~McMullen,~\emph{Complex earthquakes and Teichm\"uller
    theory}. J.Amer. Math. Soc., \textbf{11} (1998), 283--320.
\bibitem{Mess} G.~Mess,~\emph{Lorentz spacetimes of constant curvature}.
  Preprint (1990).
\bibitem{Penner}R.~C.~Penner and J.~L.~Harer,~\emph{Combiatorics of train
    tracks}. Annals of Mathematics Studies 125. Princeton University Press,
    Princeton NJ, 1992.
\bibitem{Poin} V.~Poin\'earu, A.~Fathi, F.~Laudenbach,~\emph{Travaux de Thurston
    sur le surfaces}. S\'eminaire Orsay. Ast\'erisque 66-67. Soci\'et\'e
    Math\'ematique de France, Paris, 1979.
\bibitem{Rees}M.~Rees,~\emph{An alternative approach to the ergodic theory of
    measured foliations}. Ergod.~Theory~and~Dynam.~Sys., \textbf{1} (1981), 461--488.
\bibitem{Scannell} K.~Scanell,~\emph{Flat conformal structures and the
    classification of de Sitter manifolds}. Comm. Anal. Geom., \textbf{7}
    (1999), 325--345.
\bibitem{Thurston} W.~Thurston,~\emph{Geometry and topology of three
    manifolds}. Lecture Notes, Princeton University, 1979.
\bibitem{Thurston:earth} W.~Thurston,~\emph{Earthquakes in two-dimensional
    hyperbolic geometry}. In \emph{Low-dimensional topology and Kleinian groups
    (Coventry/Durham, 1984)} 91--112, London Math. Lecture Notes Ser. 112,
    Cambridge University Press (1986).  
\end{thebibliography}
\end{document}